\documentstyle[12pt,amsmath]{article}
\begin{document}
\title{\bf Global solutions and Relaxation Limit to the Cauchy Problem
  of a Hydrodynamic Model for Semiconductors
 }
\author{\noindent
Yun-guang  Lu
\thanks{ the corresponding author: ylu2005@ustc.edu.cn}\\
K.K.Chen Institute for Advanced Studies  \\
Hangzhou Normal University, P. R. CHINA
}
\date{ }
\vskip 0.5cm
\baselineskip 18pt
\newtheorem{define}{{\bf Definition}}
\newcommand{\U}[1]{#1}
\newcommand{\Up}[0]{\Upsilon}
 \newcommand{\F}[2]{\frac{#1}{#2}}
 \newcommand{\DF}[2]{\dfrac{#1}{#2}}
\newcommand{\E}[0]{\varepsilon}
\newcommand{\lam}[0]{\lambda}
 \newcommand{\BE}[1]{\begin{array}{#1}}
\newcommand{\EE}[0]{\end{array}}
\makeatletter
\@addtoreset{equation}{section}
\makeatother\newcommand{\DIS}[0]{\displaystyle}
\newcommand{\Pa}[0]{\partial}
\newcommand{\G}[0]{\bigtriangledown}
\renewcommand{\theequation}{\arabic{section}.\arabic{equation}}
\newcommand{\usection}[1]{section}
\newtheorem{theorem}{\underline{Theorem}}
\newtheorem{lemma}[theorem]{\underline{Lemma}}
\newlength{\picwidth} \newlength{\picwidthadj}
\newcommand{\unipic}[2]{\vspace{#2}
\newcommand \epf{\rule{0.2cm}{0.2cm}\hfill \bigskip \par}
\setlength{\picwidth}{\the\textwidth}
\setlength{\picwidthadj}{0pt}
\addtolength{\picwidth}{\picwidthadj}
\special{dvitops:import #1 \the\picwidth #2}
\bigskip }
\maketitle
\begin{abstract} It is well-known that due to the lack of a technique to obtain the
a-priori $L^{\infty}$ estimate of the artificial viscosity solutions of the Cauchy
problem for the one-dimensional Euler-Poisson (or hydrodynamic) model for semiconductors,
where the energy equation is replaced by a pressure-density relation,
over the past three decades, all solutions of this model were obtained by using the
Lax-Friedrichs, Godounov schemes and Glimm scheme for both the initial-boundary value problem \cite{Zh1,Li} and the Cauchy problem \cite{MN1,PRV,HLY}; or by using the vanishing artificial viscosity method for the initial-boundary value problem \cite{Jo,HLYY}. In this paper, the existence of global entropy solutions, for the Cauchy problem of this model, is proved by using the vanishing artificial viscosity method. First, a special flux approximate is introduced to ensure the uniform boundedness of the electric field $E$ and the a-priori $L^{\infty}$ estimate, $ 0 < 2 \delta \leq \rho^{\E,\delta} \leq M(t), |u^{\E,\delta}| \leq M(t)$, where $M(t)$ could tend to infinity as the time $t$ tends to infinity, on the viscosity-flux approximate solutions $(\rho^{\E,\delta},u^{\E,\delta})$; Second, the compensated compactness theory is applied to prove the pointwise convergence of $(\rho^{\E,\delta},u^{\E,\delta})$ as $\E,\delta$ go to zero, and that the limit $(\rho(x,t),u(x,t)$ is a global entropy solution; Third, a technique, to apply the maximum principle to the combination of the Riemann invariants and $ \int_{-\infty}^{x} \rho^{\E,\delta}(x,t)- 2 \delta dx$, deduces the uniform $L^{\infty}$ estimate,  $ 0 < 2 \delta \leq \rho^{\E,\delta} \leq M, |u^{\E,\delta}| \leq M$, independent of the time $t$ and $\E,\delta$; Finally, as a by-product, the known compactness framework \cite{MN2,JR} is applied to show the relaxation limit, as the relation time $ \tau$ and  $\E,\delta$ go to zero,
for general pressure $P(\rho)$.
\end{abstract}
 { Key Words:  Entropy solution;  Viscosity method; Cauchy problem; Maximum principle; Hydrodynamic model for semiconductors; Relaxation limit
  \\}
\section{ \bf Introduction}
In this paper, we study the global generalized solutions and the relaxation limit of the one-dimensional isentropic Euler-Poisson model for semiconductor devices:
\begin{equation}\left\{
\label{1.1}
\begin{array}{l}
\rho_{t}+( \rho u)_{x}=0,         \\\\
       ( \rho u)_t+( \rho u^2+ P(\rho))_x= \rho E-a(x) \F{\rho u}{\tau}, \\\\
       E_x= \rho-b(x)
\end{array}\right.
\end{equation}
in the region $ (-\infty,+ \infty) \times [0,T]$,
with bounded initial data
\begin{equation}
\label{1.2} (\rho,u)|_{t=0}=(\rho_{0}(x),u_{0}(x)), \quad \lim_{|x| \to \infty} (\rho_{0}(x), u_{0}(x)) =(0,0), \quad \rho_{0}(x) \ge 0
\end{equation}
and a condition at $- \infty$  for the electric field
\begin{equation}
\label{1.3} \lim_{x \rightarrow - \infty} E(x,t)= E_{-}, \hspace{0.3cm} \mbox{ for a.e.} \quad t \in (0,\infty),
\end{equation}
where $T, E_{-}$ are fixed constants,  $\rho \geq 0$ denotes the electron density, $u$ the (average) particle velocity and $E$ the electric field, which is generated by the Coulomb force of the
particles. The two given functions $a(x), b(x)$, respectively represent a damping coefficient and
the concentration of a fixed background charge \cite{PRV}. The pressure-density relation is $P=P(\rho)= \F{1}{ \gamma} \rho^{\gamma} $, where $ \gamma \geq 1$ corresponds to the adiabatic exponent
and $ \tau >0 $  is the momentum relaxation time.

As a simplified hydrodynamic model, system (\ref{1.1}) was first derived in \cite{DM}. The existence
of a unique smooth solution from the steady-state of (\ref{1.1}) in the subsonic case was proved in
\cite{DM}. The existence of a local smooth solution of the time-dependent problem (\ref{1.1}) was
proved by using Lagrangian mass coordinates in \cite{Zh2}.  Due to the formation of
shocks, one cannot expect to obtain global smooth solution in the general case. For
example, the numerical simulation of a steady-state shock wave in the hydrodynamic model was first
presented by Gardner in \cite{Ga}.  About the existence of the global weak solutions of the time-dependent problem
(\ref{1.1}), all works are concentrated in the two respects. First, the solutions were obtained by using the Lax-Friedrichs and Godounov schemes, or Glimm scheme for both the initial-boundary value problem \cite{Zh1,Li} and the Cauchy problem \cite{MN1,PRV,HLY}; Second, the solutions were obtained by using the vanishing artificial viscosity method for the initial-boundary value problem \cite{Jo,HLYY}. All these results  are based on the corresponding compactness framework on the following homogeneous system of isentropic gas dynamics
\begin{equation}\left\{
\label{1.4}
\begin{array}{l}
\rho_{t}+( \rho u)_{x}=0,         \\\\
       ( \rho u)_t+( \rho u^2+ P(\rho))_x=0.
\end{array}\right.
\end{equation}
More interesting existence or non-existence results  on related hydrodynamic model of semiconductor devices can be found in \cite{TW,Pe1,Xu,PX,FXZ,GN,GS,HPY,HZ,LY,LNX,Ts,Wang,WY} and the references cited therein. How to obtain a global solution of the
Cauchy problem of system (\ref{1.1}) by using the vanishing artificial viscosity method is an open problem
in the last three decades. The main difficulties are in the following.

The classical vanishing viscosity method is to add the viscosity terms
to the right-hand side of system (\ref{1.1}) and consider the
problem for the related system
\begin{equation}\left\{
\label{1.5}
\begin{array}{l}
 \rho_{t}+( \rho u)_{x} =\E \rho_{xx},         \\\\
 ( \rho u)_t+( \rho u^2+ P(\rho))_x= \E
 (\rho u)_{xx}+ \rho E-  a(x) \F{\rho u}{\tau}, \\\\
 E_x= \rho-b(x).
  \end{array}\right.
\end{equation}
If we consider $(\rho,m)$, where $m=\rho u$ as two
independent variables, then the
term $ \rho u^2 = \F{m^2}{\rho}$ in (\ref{1.5}) is singular near the line $\rho=0$. By using the energy method given in
\cite{Di1}, the maximum priciple in \cite{Lu4} or with the help of the Green function in \cite{Lu3,Pe2}, we may obtain the
positive lower bound $  \rho^{\E} \geq  c(t,c_0,\E) > 0$, if $ \rho_{0}(x) \geq c_0 > 0$,
where $c_0 $ is a positive constant and $c(t,c_0,\E)$ could tend to
zero as the time $t$ tends to infinity or $\E$ tends to zero. However, the lower bound $ c(t,c_0,\E)$ could not
ensure the boundedness of the electric field $ E^{\E} $ because  $ E^{\E}(x,t)= E_{-}+ \int_{- \infty}^{x} \rho(x,t) dx
- \int_{- \infty}^{x} b(x) dx$ from the third equation in (\ref{1.5}).

To overcome the above difficulties, in this paper, we apply the combination of the flux approximate coupled with the
classical vanishing viscosity \cite{Lu1,Lu2} to study the following system
\begin{equation}\left\{
\label{1.6}
\begin{array}{l}
 \rho_{t}+( (\rho-2 \delta) u)_{x} =\E \rho_{xx},         \\\\
 ( \rho u)_t+( \rho u^2- \delta u^{2}+ P_{1}( \rho, \delta))_x= \E
 (\rho u)_{xx}+ \rho E- a(x) \F{\rho u}{\tau}, \\\\
 E_x= (\rho- 2 \delta) -b(x)
  \end{array}\right.
\end{equation}
with the initial data
\begin{equation}
\label{1.7}
      (\rho^{\E,\delta}(x,0),u^{\E,\delta}(x,0))=(\rho_{0}(x)+ 2
      \delta,u_{0}(x)) \ast G^{\E},
\end{equation}
where $(\rho_{0}(x),u_{0}(x))$ are given in (\ref{1.2}), $ \delta > 0$  denotes a regular perturbation constant,
the perturbation pressure
 \begin{equation}
\label{1.8}
P_{1}( \rho, \delta)= \int_{2 \delta}^{\rho} \F{t-2
\delta}{t}P'(t)dt,
\end{equation}
$G^{\E}$ is a mollifier such that $(\rho^{\E,\delta}(x,0),u^{\E,\delta}(x,0)) $
are smooth and
\begin{equation}
\label{1.9}
   \lim_{|x| \to \infty} (\rho^{\E,\delta}(x,0),u^{\E,\delta}(x,0)) =( 2 \delta,0),
   \quad \lim_{|x| \to \infty} (\rho^{\E,\delta}_{x}(x,0),u^{\E,\delta}_{x}(x,0)) =(0,0).
\end{equation}
By simple calculations, two eigenvalues of system (\ref{1.6}) are
\begin{equation}
\label{1.10} \lam^{\delta}_{1}= \F{m}{ \rho}- \F{\rho-2
\delta}{\rho} \sqrt{P'( \rho)}, \quad \lam^{\delta}_{2}= \F{m}{
\rho}+ \F{\rho-2 \delta}{\rho} \sqrt{P'( \rho)}
\end{equation}
with corresponding two Riemann invariants
\begin{equation}
\label{1.11} z(\rho,u)=  \int_{l}^{ \rho} \F{ \sqrt{P'(s)}}{s}ds -u,
\quad w(\rho,u)= \int_{l}^{ \rho} \F{ \sqrt{P'(s)}}{s}ds+u,
\end{equation}
where $l$ is a constant.

By using the first equation in (\ref{1.6}), we have $  \rho^{\E,\delta}(x,t) \geq 2  \delta $
and
\begin{equation}
\label{1.12}
\int_{- \infty}^{\infty} \rho^{\E,\delta}(x,t) - 2  \delta dx \leq
\int_{- \infty}^{\infty} \rho^{\E}_{0}(x) - 2  \delta dx =
\int_{- \infty}^{\infty} \rho_{0}(x) dx \leq M;
\end{equation}
by using the third equation in (\ref{1.6}),
\begin{equation}
\label{1.13}
|E^{\E,\delta}(x,t)| =|E_{-}+ \int_{- \infty}^{x} \rho^{\E,\delta}(x,t) - 2  \delta dx
- \int_{- \infty}^{x} b(x) dx| \leq M,
\end{equation}
where $ M $ denotes a suitable positive constant.

With the help of the estimate (\ref{1.13}), we immediately have the first existence result in this paper
\begin{theorem} Let the initial data $( \rho_{0}(x),u_{0}(x)) $ be bounded measurable, $\rho_0(x) \geq 0,
|\rho_{0}(x)|_{L^{1}(R)} \leq M$ and $ a(x) \in C^{2}(R), 0 \leq a(x) \leq M, |b(x)|_{L^{1}(R)} \leq M$. Let
$ P(\rho)\in C^{2}(0,\infty), P'(\rho)>0,  P''(\rho) > 0$
for $\rho>0$; and
\begin{equation}
\label{1.14}
\int_{c}^{\infty}\F{\sqrt{P'(\rho)}}{\rho}d\rho=\infty, \quad
\int_{0}^{c}\F{\sqrt{P'(\rho)}}{\rho}d\rho<\infty,
\quad \forall c>0.
\end{equation}
Then for fixed $ \E >0, \delta >0 $, the smooth solution $(\rho^{\E, \delta}(x,t), u^{\E,\delta}(x,t),E^{\E,\delta}(x,t))$ of the  problem (\ref{1.3}),
(\ref{1.6}) and  (\ref{1.7})
exists in any region $(- \infty, \infty)  \times [0,T), T >0$,  satisfies (\ref{1.12}), (\ref{1.13}),
\begin{equation}
\label{1.15}
  \lim_{|x| \to \infty} (\rho^{\E,\delta}(x,t), u^{\E,\delta}(x,t)) =( 2 \delta,0), \quad
  \lim_{|x| \to \infty} (\rho_{x}^{\E,\delta}(x,t), m_{x}^{\E,\delta}(x,t)) =(0,0)
  \end{equation}
  and
\begin{equation}
\label{1.16}
0 < 2 \delta \leq  \rho^{\E,\delta}(x,t) \leq M(T), \quad
 |u^{\E,\delta}(x,t)|= |\F{m^{\E,\delta}(x,t)}{ \rho^{\E,\delta}(x,t)}| \leq M(T),
  \end{equation}
where $M(T)$ is a positive constant, being independent
of $\tau, \E$ and $ \delta$, but could tend to infinity
as $T$ tends to infinity.

Particularly, for the polytropic gas and $ \gamma >1$, there exists a subsequence
(still labelled) $(\rho^{\E,\delta}(x,t), \rho^{\E,\delta}(x,t) u^{\E,\delta}(x,t),E^{\E,\delta}(x,t))$
which converges almost everywhere on any bounded and open set
$ \Omega \subset R \times R^+$:
\begin{equation}
\label{1.17}
 (\rho^{\E,\delta}(x,t), \rho^{\E,\delta}(x,t) u^{\E,\delta}(x,t),E^{\E,\delta}(x,t))
  \rightarrow ( \rho(x,t),\rho(x,t) u(x,t), E(x,t)),
\end{equation}
 as  $ \E \downarrow 0^{+}, \delta \downarrow 0^{+}$, where the limit
$(\rho(x,t),\rho(x,t)u(x,t),E(x,t))$ is a weak solution
 of the  problem (\ref{1.1})-(\ref{1.3}).

Furthermore, for the isothermal gas  $ \gamma = 1 $, for fixed $ \E >0, \delta >0 $, the smooth viscosity solution $(\rho^{\E, \delta}(x,t), u^{\E,\delta}(x,t),E^{\E,\delta}(x,t))$, of the problem (\ref{1.3}),
(\ref{1.6}) and  (\ref{1.7}),
exists in any region $(- \infty, \infty)  \times [0,T), T >0$,  satisfies (\ref{1.12}), (\ref{1.13}),
\begin{equation}
\label{1.18}
  \lim_{|x| \to \infty} (\rho^{\E,\delta}(x,t), u^{\E,\delta}(x,t)) =( 2 \delta,0), \quad
  \lim_{|x| \to \infty} (\rho_{x}^{\E,\delta}(x,t), m_{x}^{\E,\delta}(x,t)) =(0,0)
  \end{equation}
  and
  \begin{equation}
\label{1.19}
0 < 2 \delta \leq  \rho^{\E,\delta}(x,t) \leq M(T), \quad
 |\rho^{\E,\delta}(x,t) u^{\E,\delta}(x,t)|= |m^{\E,\delta}(x,t)| \leq M(T).
  \end{equation}
Moreover,  there also exists a subsequence
of $(\rho^{\E,\delta}(x,t), \rho^{\E,\delta}(x,t) u^{\E,\delta}(x,t),E^{\E,\delta}(x,t))$,
whose limit
$(\rho(x,t),\rho(x,t)u(x,t),E(x,t))$ is a weak solution of the problem (\ref{1.1})-(\ref{1.3}).
\end{theorem}
\begin{define} $(\rho(x,t),u(x,t),E(x,t))$  is called a weak entropy
solution of the problem (\ref{1.1})-(\ref{1.3}) if
\begin{equation}\left\{
\label{1.20}
\begin{array}{l}
      \int_{0}^{\infty} \int_{- \infty}^{\infty} \rho \phi_{t}+( \rho u) \phi_{x} dxdt
      + \int_{- \infty}^{\infty} \rho_{0}(x) \phi(x,0) dx=0,  \\  \\
       \int_{0}^{\infty} \int_{- \infty}^{\infty}  \rho u \phi_t+( \rho u^2+ P( \rho)) \phi_x
     + ( \rho E -a(x) \F{\rho u}{\tau} ) \phi dxdt \\\\
      + \int_{- \infty}^{\infty} \rho_{0}(x) u_{0}(x) \phi(x,0) dx=0,  \\\\
      \int_{0}^{\infty} \int_{- \infty}^{\infty} E \phi_{x}+ (\rho-b(x)) \phi dxdt
\end{array}\right.
\end{equation}
holds for all test function $ \phi \in C_{0}^{1}(R \times R^{+})$
and
\begin{equation}\begin{array}{ll}
\label{1.21} \int_0^{ \infty} \int_{ -\infty}^{\infty} \eta(\rho,m)
\phi_t+ q( \rho,m) \phi_x + ( \rho E -a(x) \F{\rho u}{\tau} ) \eta(\rho,m)_{m}  \phi dxdt \geq 0
\end{array}
\end{equation}
holds for any non-negative test function $ \phi \in C_0^{ \infty}(R
\times R^{+}- \{t=0\}),$ where  $m= \rho u $ and $ ( \eta, q)$ is a
pair of convex entropy-entropy flux of system (\ref{1.4}).
\end{define}
When we study the  limit \cite{MN2,JR,JP} of $(\rho(x,t),\rho(x,t)u(x,t),E(x,t))$ as the relaxation time $\tau \rightarrow 0$ or similarly consider the large time behavior \cite{HPY,HLYY,Yu} of $(\rho(x,t),\rho(x,t)u(x,t),E(x,t))$  as the time $t \rightarrow \infty $, we need to prove that the upper bounds, given in (\ref{1.16}) and (\ref{1.19}), are independent of $t$.

Based on the uniform bound assumption on the approximate solutions constructed by using the
fractional step Lax-Friedrichs scheme and Godounov scheme, the authors in \cite{MN2} proved that the limit
  $ (N(x,s), J(x,s), \Up (x,s))$ of the sequence $(N^{\tau}(x,s),J^{\tau}(x,s), \Upsilon^{\tau}(x,s))$, as $\tau \downarrow 0^{+}$,
 is a solution of the following well-known drift-diffusion equations
 \begin{equation}\left\{
\label{1.22}
\begin{array}{l}
      N_{s}+J_{x}=0      \\\\
       P(N)_{x}= N \Upsilon - a(x)J, \\\\
\Upsilon_x=N-b(x)
\end{array}\right.
\end{equation}
in the sense of distributions, where $ (N^{\tau}(x,s),J^{\tau}(x,s), \Upsilon^{\tau}(x,s))$ are obtained by introducing the scaled variables on the solution $(\rho(x,t),\rho(x,t)u(x,t),E(x,t))$  of the problem (\ref{1.1})-(\ref{1.3}),
\begin{equation}
\label{1.23}
N^{\tau}(x,s)= \rho(x,\F{s}{\tau}),  \quad J^{\tau}(x,s)=\F{1}{\tau} m(x,\F{s}{\tau}),
\quad \Up^{\tau}(x,s)= E(x,\F{s}{\tau}).
  \end{equation}
After giving up the attempt to obtain the uniform bound on the approximate solutions, the authors in \cite{JP}
 constructed a family of positive and convex entropies
to deduce the high energy estimates of solutions and the uniform $L^{p}, 1 \leq p < \infty$, estimates of
 the approximate solutions. Based on the $L^{p}$ estimates and the technical assumption $\gamma= 1+ \F{2}{m}, m \geq 1$ being an
 integer, the zero relaxation limit of $ (N^{\tau}(x,s),J^{\tau}(x,s), \Upsilon^{\tau}(x,s))$ was proved by using the
 compensated compactness method.

For the isothermal gas $ \gamma=1,$  under the assumptions $u_{0}(x) \in BV (R)$ and $\ln \rho_{0}(x) \in BV (R)$,
where $ BV (R)$ is the space of functions with bounded variation, the relaxation limit was proved in \cite{JR} by introduing
the Glimm scheme \cite{Gl} to construct the approximate solutions of the problem (\ref{1.1})-(\ref{1.3}).

The second purpose of this paper is to prove the uniform $ L^{\infty}$ estimates, of
 the flux-viscosity solutions $(\rho^{\E,\delta}(x,t), \rho^{\E,\delta}(x,t) u^{\E,\delta}(x,t),E^{\E,\delta}(x,t))$ of the problem (\ref{1.3}), (\ref{1.6}) and (\ref{1.7}), independent of $\E, \delta, \tau $ and the time $t$,  and to prove that
 the limit  $ (N(x,s), J(x,s), \Up (x,s))$ of
 \begin{equation}
\label{1.24}
N^{\tau}(x,s)= \rho^{\E,\delta}(x,\F{s}{\tau}),  \quad J^{\tau}(x,s)=\F{1}{\tau} m^{\E,\delta}(x,\F{s}{\tau}),
\quad \Up^{\tau}(x,s)= E^{\E,\delta}(x,\F{s}{\tau})
  \end{equation}
  as the parameters $\E,\delta,\tau$ go to zero, is a generalized solution of the drift-diffusion equations (\ref{1.22}) in the sense of distributions. Precisely, we have the following Theorems 2 and 3.
  \begin{theorem} Let the conditions on the initial data, $a(x)$, $b(x)$  in Theorem 1 be satisfied. Suppose  $ b(x) \geq 0$,
 \begin{equation}
\label{1.25}
\int_{- \infty}^{\infty} b(x) dx < E_{-}
\end{equation}
and $a(x), C(x)$  satisfy
\begin{equation}
\label{1.26} 0 < a_{0} \leq a(x) \leq M, \quad a'(x) \leq 0, \quad C'(x) \geq 0,
\end{equation}
where $a_{0} $  is a positive constant,
\begin{equation}
\label{1.27}
C(x)=\F{1}{a(x)}(E_{-}- \int_{- \infty}^{x} b(x) dx).
\end{equation}
Then, the smooth viscosity solution $(\rho^{\E, \delta}(x,t), u^{\E,\delta}(x,t),E^{\E,\delta}(x,t))$ of the problem
(\ref{1.3}), (\ref{1.6}) and (\ref{1.7}) satisfies (\ref{1.13}) and
\begin{equation}
\label{1.28}
0 < 2 \delta \leq  \rho^{\E,\delta}(x,t) \leq M, \quad
 |u^{\E,\delta}(x,t)|= |\F{m^{\E,\delta}(x,t)}{ \rho^{\E,\delta}(x,t)}| \leq M
  \end{equation}
if $P(\rho)= \rho^{\gamma}, \gamma > 1$, where $M$ is a positive constant, being independent
of $\E, \delta, \tau$ and the time $T$; and
\begin{equation}
\label{1.29}
0 < 2 \delta \leq  \rho^{\E,\delta}(x,t) \leq e^{2 M}, \mbox{ }
  \ln \rho^{\E,\delta}(x,t) -M \leq u^{\E,\delta}(x,t) \leq M- \ln \rho^{\E,\delta}(x,t)
  \end{equation}
if $P(\rho)= \rho$.
\end{theorem}
\begin{theorem} Suppose the uniform estimates (\ref{1.28}) and (\ref{1.29}) in Theorem 2 are true. Let
\begin{equation}
\label{1.30}
      \eta^{\star}(\rho,m)= \F{m^2}{2 \rho} + \rho \int_{2 \delta}^{\rho} \F{P(s)}{s^2}ds
\end{equation}
and assume
\begin{equation}
\label{1.31}
	\int^{\infty}_{-\infty}  \eta^{\star}(\rho^{\E, \delta}_{0}(x), u^{\E, \delta}_{0}(x)) dx \leq M.
\end{equation}
Then there exists a subsequence
(still labelled) $(N^{\tau}(x,s), J^{\tau}(x,s),\Up^{\tau}(x,s))$
which converges almost everywhere on any bounded and open set
$ \Omega \subset R \times R^+$:
\begin{equation}
\label{1.32}
 (N^{\tau}(x,s), J^{\tau}(x,s),\Up^{\tau}(x,s))
  \rightarrow ( N(x,s),J(x,s), \Up(x,s)),
\end{equation}
 as $ \E = o( \sqrt{P'( 2 \delta)} ) $ and  $ \tau \downarrow 0^{+}, \delta \downarrow 0^{+}$, where the limit
$( N(x,s),J(x,s), \Up(x,s))$ is a bounded weak solution
 of the drift-diffusion equations (\ref{1.22}) in the sense of distributions, and $ \Up(x,s) $ is Lipschitz continuous with respect to the space variable $x$.
\end{theorem}

{\bf Remark 1.} From the proof of Theorem 2 in Section 3, we can see that, if $b(x) \leq 0$,  the condition (\ref{1.25}) can be removed.
A simple example, to ensure (\ref{1.26}) and (\ref{1.27}),  is to let $a(x)= E_{-}- a_{1}
\int_{- \infty}^{x} b(x) dx$, where $ a_{1} \geq 1 $ is a constant. In general $C(x)$ satisfies
\begin{equation}
\label{1.33} 0 < \F{1}{M}(E_{-}- \int_{- \infty}^{\infty} b(x) dx) \leq
C(x) \leq \F{1}{a_{0}}E_{-}.
\end{equation}

\section{ \bf Proof of Theorem 1}
To prove Theorem 1, we first have the following local existence lemma.
\begin{lemma} Let the conditions in Theorem 1 be satisfied. Then: {\it (1)} for any fixed $ \E, \delta, \tau >0$, the problem (\ref{1.3}),(\ref{1.6}) and (\ref{1.7})
 always has a local smooth solution $(\rho^{\E,\delta}, u^{\E,\delta}) \in C^{ \infty} (R \times (0, t_{1}))$
for a small time $ t_{1}$, which depends only on the $L^{ \infty} $ norm
of the initial data $(\rho_{0}(x), u_{0}(x))$; and  {\it (2)} the solution $(\rho^{\E,\delta}, u^{\E,\delta})$ satisfies:
\begin{equation}
\label{2.1}
  \lim_{|x| \to \infty} (\rho^{\E,\delta}(x,t), u^{\E,\delta}(x,t)) =( 2 \delta,0), \quad
  \lim_{|x| \to \infty} (\rho_{x}^{\E,\delta}(x,t), m_{x}^{\E,\delta}(x,t)) =(0,0)
  \end{equation}
 and
\begin{equation}
\label{2.2}
\delta \leq \rho^{\E,\delta} \leq 2 M, \quad |\rho^{\E,\delta} u^{\E,\delta}| \leq  2 M,
 \quad \int_{- \infty}^{\infty} \rho^{\E,\delta}(x,t) - 2  \delta dx \leq M.
\end{equation}
  \end{lemma}
{ \bf Proof of Lemma 4.}  First, we integrate the third equation in (\ref{1.6}) over $(-\infty, x)$  and replace
$E^{\E,\delta}$ in the second equation by
\begin{equation}
\label{2.3}
E_{-}+ \int_{- \infty}^{x} \rho^{\E,\delta}(x,t) - 2  \delta dx - \int_{- \infty}^{x} b(x) dx.
\end{equation}
Then, the local existence result in {\em (1)} can be
obtained by applying the contraction mapping principle
to an integral representation for a solution, following the standard
theory of semilinear parabolic systems. In fact, by applying the Green function we may obtain the following integral
 representation of the first two equations in (\ref{1.6}),
\begin{equation}\left\{
\label{2.4}
\begin{array}{l}
       \rho(x,t)= \rho^{0}(x,t)+ \int_{0}^{t} \int_{-\infty}^\infty (\rho(y,s)- 2 \delta)u(y,s) G_{y}(x-y,t-s) dy ds,       \\\\

      \rho(x,t) u(x,t)= m^{0}(x,t)+ \int_{0}^{t} \int_{-\infty}^\infty f(\rho(y,s),u(y,s)) G_{y}(x-y,t-s) \\\\
      + (\rho (\int_{- \infty}^{x} \rho(x,t) - 2  \delta dx - \int_{- \infty}^{x} b(x) dx) - a(x) \F{\rho u}{\tau})
       G(x-y,t-s) dy ds,
\end{array}\right.
\end{equation}
where
\begin{equation}\left\{
\label{2.5}
\begin{array}{l}
\rho^{0}(x,t)= \int_{-\infty}^\infty \rho^{\E,\delta}_{0}(y) G(x-y,t) dy,  \quad m^{0}(x,t)= \int_{-\infty}^\infty \rho^{\E,\delta}_{0}(y) u^{\E,\delta}_{0}(y) G(x-y,t) dy, \\\\
f(\rho,u) =  \rho u^2- \delta u^{2}+ P_{1}( \rho, \delta), \quad G(x,t)=\F{1}{\sqrt{\pi\E t}}exp(-\F{x^{2}}{4\E t}).
\end{array}\right.
\end{equation}
Second, we construct the iteration sequence $(\rho^{(n)}(x,t), u^{(n)}(x,t)), n \geq1, $ as follows
\begin{equation}\left\{
\label{2.6}
\begin{array}{l}
       \rho^{(n)}(x,t)= \rho^{0}(x,t)+ \int_{0}^{t} \int_{-\infty}^\infty (\rho^{(n-1)}(y,s)- 2 \delta)u^{(n-1)}(y,s) G_{y}(x-y,t-s) dy ds,       \\\\

      \rho^{(n)}(x,t) u^{(n)}(x,t)= m^{0}(x,t)+ \int_{0}^{t} \int_{-\infty}^\infty f(\rho^{(n-1)}(y,s),u^{(n-1)}(y,s)) G_{y}(x-y,t-s) \\\\
      + (\rho^{(n-1)} (\int_{- \infty}^{x} \rho^{(n-1)}(x,t) - 2  \delta dx - \int_{- \infty}^{x} b(x) dx)\\\\
       - a(x) \F{\rho^{(n-1)} u^{(n-1)}}{\tau})
       G(x-y,t-s) dy ds,
\end{array}\right.
\end{equation}
where, when $n-1=0,  \rho^{(0)}(x,t)= \rho^{\E,\delta}_{0}(x,0), u^{(0)}(x,t)=u^{\E,\delta}_{0}(x,0)$.

Third, since the initial data satisfy (\ref{1.9}) and
$ 0 < 2 \delta \leq \rho^{\E,\delta}_{0}(x,0) \leq M, |m^{\E,\delta}_{0}(x,0)|=| \rho^{\E,\delta}_{0}(x,0) u^{\E,\delta}_{0}(x,0) | \leq M$, by induction, we can easily prove that there exists a small $t_{1}$ such that the following estimates are
true for any $ n \geq 1$ and for $ t \in [0, t_{1}]$,
\begin{equation}\left\{
\label{2.7}
\begin{array}{l}
\delta  \leq \rho^{(n)}(x,t) \leq 2 M, \quad |m^{(n)}(x,t)|=| \rho^{(n)}(x,t) u^{(n)}(x,t)| \leq 2 M, \\\\
\sup\limits_{ 0 \leq t \leq t_{1} } (|\rho^{(n)}(\cdot,t)-\rho^{(n-1)}(\cdot,t)|_{L^{\infty}(R)}+|m^{(n)}(\cdot,t)-m^{(n-1)}(\cdot,t)|_{L^{\infty}(R)})
\\\\ \leq c \sup\limits_{ 0 \leq t \leq t_{1} } (|\rho^{(n-1)}(\cdot,t)-\rho^{(n-2)}(\cdot,t)|_{L^{\infty}(R)}+|m^{(n-1)}(\cdot,t)-m^{(n-2)}(\cdot,t)|_{L^{\infty}(R)}),
  \\\\
\lim_{|x| \to \infty} (\rho^{(n)}(x,t),m^{(n)}(x,t)) =( 2 \delta,0),
   \quad \lim_{|x| \to \infty} (\rho^{n}_{x}(x,t),m^{n}_{x}(x,t)) =(0,0),
\end{array}\right.
\end{equation}
where $ c$ is a constant, $0 < c < 1$.

Furthermore, if $ |\rho^{(n-1)}(\cdot,t) - 2 \delta|_{L^{1}(R)} \leq M$, we have
\begin{equation}\begin{array}{ll}
\label{2.8}
|\rho^{(n)}(\cdot,t) - 2 \delta|_{L^{1}(R)} \leq | \rho^{0}(x,t) - 2 \delta|_{L^{1}(R)} \\\\
+ |\int_{0}^{t} \int_{-\infty}^\infty ((\rho^{(n-1)}(y,s)- 2 \delta)u^{(n-1)}(y,s))_{y}
 (\int_{-\infty}^\infty G(x-y,t-s) dx) dy ds|  \\\\
 \leq M + |\int_{0}^{t} \int_{-\infty}^\infty ((\rho^{(n-1)}(y,s)- 2 \delta)u^{(n-1)}(y,s))_{y}
  dy ds|= M
\end{array}
\end{equation}
due to
\begin{equation}
\label{2.9}
\int_{-\infty}^\infty G(x-y,t-s) dx=1, \quad \lim_{|x| \to \infty} (\rho^{(n-1)}(x,t),m^{(n-1)}(x,t)) =( 2 \delta,0).
\end{equation}
Thus there exists a solution $(\rho(x,t),m(x,t))$, for the integral equations in (\ref{2.4}), which satisfies (\ref{2.2}) and
\begin{equation}
\label{2.10}
  \lim_{|x| \to \infty} (\rho^{\E,\delta}(x,t), m^{\E,\delta}(x,t)) =( 2 \delta,0), \quad
  \lim_{|x| \to \infty} (\rho_{x}^{\E,\delta}(x,t), m_{x}^{\E,\delta}(x,t)) =(0,0).
  \end{equation}
Since $ \lim_{|x| \to \infty} u(x,t) = \lim_{|x| \to \infty} \F{m(x,t)}{\rho}=0$,
we obtain (\ref{2.1}) and so {\bf the proof of Lemma 4}.
More details about the local solution of a given parabolic system  can be seen in \cite{LSU, Sm} or \cite{Lu5}.

Whenever we have an {\it a priori} $L^{ \infty} $ estimate of the
local solution, it is clear that the local time $ t_{1}$
can be extended to $T$ step by step since the step time depends only on the
 $L^{ \infty} $ norm.

 { \bf Proof of Theorem 1.} By using the maximum principle  to the first equation in (\ref{1.6}), we have
 $  \rho^{\E,\delta}(x,t) \geq 2  \delta $ (See Lemma 2.2 in \cite{Lu6} for the details),  and
\begin{equation}
\label{2.11}
\int_{- \infty}^{\infty} \rho^{\E,\delta}(x,t) - 2  \delta dx \leq
\int_{- \infty}^{\infty} \rho^{\E,\delta}_{0}(x) - 2  \delta dx =
\int_{- \infty}^{\infty} \rho_{0}(x) dx \leq M
\end{equation}
(See Theorem 1 in \cite{Lu7} for the details).

By using the third equation in (\ref{1.6}), we have (\ref{2.3}) and so
\begin{equation}
\label{2.12}
|E^{\E,\delta}(x,t)| =|E_{-}+ \int_{- \infty}^{x} \rho^{\E,\delta}(x,t) - 2  \delta dx
- \int_{- \infty}^{x} b(x) dx| \leq M_{1}
\end{equation}
for a suitable constant $M_{1}>0$.

We  multiply the first two equations in  (\ref{1.6}) by  $(z_{\rho},
z_{m})$ and $(w_{\rho}, w_{m})$, respectively, where $(z,w)$ are given in (\ref{1.11}), to obtain
\begin{equation}\begin{array}{ll}
\label{2.13}
 z_{t}+ \lam^{\delta}_{1}z_{x}= -E+ a(x) \F{u}{\tau} \\\\
 +  \E z_{xx}+ \F{2 \E}{\rho} \rho_{x}z_{x}- \F{\E}{2 \rho^{2}
   \sqrt{P'(\rho)}}(2P'+ \rho P'') \rho_{x}^{2}
   \end{array}
\end{equation}
and
\begin{equation}\begin{array}{ll}
\label{2.14}
 w_{t}+ \lam^{\delta}_{2}w_{x}= E- a(x) \F{u}{\tau} \\\\
+\E w_{xx}+ \F{2 \E}{\rho} \rho_{x}w_{x}- \F{\E}{2 \rho^{2}
   \sqrt{P'(\rho)}}(2P'+ \rho P'') \rho_{x}^{2}.
\end{array}
\end{equation}
Since $|E^{\E,\delta}(x,t) | \leq M_{1}$, we may let
\begin{equation}
\label{2.15}
w= w_{1}+ M_{1} t,  \quad z= z_{1}+ M_{1} t
\end{equation}
to obtain from (\ref{1.11}), (\ref{2.13}) and (\ref{2.14}) that
\begin{equation}\left\{
\label{2.16}
\begin{array}{l}
z_{1t}+ \lam^{\delta}_{1}z_{1x}
\leq a(x) \F{w_{1}-z_{1}}{2 \tau}
 +  \E z_{1xx}+ \F{2 \E}{\rho} \rho_{x}z_{1x},  \\\\
 w_{1t}+ \lam^{\delta}_{2}w_{1x} \leq - a(x) \F{w_{1}-z_{1}}{2 \tau}
+\E w_{1xx}+ \F{2 \E}{\rho} \rho_{x}w_{1x}.
\end{array}\right.
\end{equation}
Since $\tau > 0, a(x) \geq 0$, we may apply the maximum principle to (\ref{2.16}) to obtain
$z_{1}(\rho^{\E,\delta},u^{\E,\delta}) \leq M_{2}, w_{1}(\rho^{\E},u^{\E}) \leq M_{2}$ if
the initial data satisfy the same estimates. Then
\begin{equation}
\label{2.17}
z(\rho^{\E,\delta},u^{\E,\delta}) \leq M_{2}+ M_{1} t, \quad
w(\rho^{\E,\delta},u^{\E,\delta}) \leq M_{2}+ M_{1} t,
\end{equation}
which give us the estimates in (\ref{1.16}) if $P(\rho)$ satisfies (\ref{1.14}), and the estimates in (\ref{1.19}) if
$P(\rho)= \rho$. Thus, for fixed $ \E >0, \delta >0 $, the smooth viscosity solution $(\rho^{\E, \delta}(x,t), u^{\E,\delta}(x,t),E^{\E,\delta}(x,t))$ of the problem (\ref{1.3}),(\ref{1.6}) and (\ref{1.7})
exists in any region $(- \infty, \infty)  \times [0,T), T >0$.

With the help of the positive lower bound $ \rho^{\E,\delta} \geq 2 \delta $ again, we proved in \cite{Lu2} that
\begin{equation}
\label{2.18}
\eta_{t}( \rho^{\E,\delta}(x,t),m^{\E,\delta}(x,t))
+q_{x}( \rho^{\E,\delta}(x,t),m^{\E,\delta}(x,t))
\end{equation}
are compact in  $H^{-1}_{loc}(R \times R^{+})$, for general pressure function $P(\rho)$, where $(\eta,q)$ is any
weak entropy-entropy flux pair of (\ref{1.4}), as
$\E$ and $\delta $ tend
to zero, with respect to the viscosity solutions
$( \rho^{\E,\delta}(x,t),m^{\E,\delta}(x,t)$
of the problem (\ref{1.3}), (\ref{1.6}) and (\ref{1.7}). Thus,  we obtain the
pointwise convergence  $( \rho^{\E,\delta}(x,t),m^{\E,\delta}(x,t)) \rightarrow (\rho(x,t),m(x,t))$ by using the compactness
frameworks given in \cite{LPS,LPT,Di1,Di2,DCL1,DCL2} when $P(\rho)= \F{1}{ \gamma} \rho^{\gamma}, \gamma >1 $ and in \cite{HW} when $\gamma=1$.

Furthermore, by  the third equation in (\ref{1.6}) and (\ref{2.12}), $E^{\E,\delta}_{x} $ are uniformly bounded in $L^{1}_{loc}(R \times R^{+}) $
 and also bounded in $ w^{-1, p}_{loc}(R \times R^{+}), p >2$,  then $E^{\E,\delta}_{x} $ are compact in $H^{-1}_{loc}(R \times R^{+})$ by using the Murat's lemma \cite{Mu}.

Moreover, by (\ref{2.12}) and the first equation in (\ref{1.6}), we have
\begin{equation}
\label{2.19}
E^{\E,\delta}_{t}= \int_{-\infty}^{x} \rho^{\E,\delta}_{t} dx= - (  \rho^{\E,\delta}- 2 \delta) + \E \rho^{\E,\delta}_{x},
\end{equation}
which are clearly compact in $H^{-1}_{loc}(R \times R^{+})$.

Thus, we may apply the Div-Curl lemma \cite{Ta} to the pairs of functions
\begin{equation}
\label{2.20}
(0,E^{\E,\delta}),    \quad (E^{\E,\delta},0),
\end{equation}
to obtain
\begin{equation}
 \label{2.21}
\overline{ E^{\E,\delta}} \cdot \overline{ E^{\E,\delta}}=  \overline{
(E^{\E,\delta})^{2}},
\end{equation}
where $ \overline{f(s^{\E,\delta})}$ denotes the weak-star limit of
$f(s^{\E,\delta})$, which deduces the pointwise convergence of $E^{\E,\delta}$,
\begin{equation}
\label{2.22}
E^{\E,\delta}(x,t) \rightarrow E(x,t),  \mbox{ strongly in }  L^{p}(\Omega),  \mbox{ for all }  p \geq 1,
\mbox{ as } \E \rightarrow 0, \delta \rightarrow 0.
\end{equation}
Letting $\E,\delta$ in (\ref{1.6}) go to zero, we may prove that the limit $(\rho,u,E)$ satisfies (\ref{1.20})-(\ref{1.21}),
and so complete the proof of Theorem 1.

\section{ \bf Proof of Theorem 2}
To obtain the bound, of $z,w$, independent of the time $t$, we rewrite
\begin{equation}\begin{array}{ll}
\label{3.1}
E(x,t) = E_{-} + \int_{- \infty}^{x} \rho(x,t) - 2  \delta dx - \int_{- \infty}^{x} b(x) dx
\\\\= a(x)[\F{1}{a(x)}
\int_{- \infty}^{x} \rho(x,t) - 2  \delta dx
+ C(x)]
= a(x) A(x,t),
\end{array}
\end{equation}
where
\begin{equation}
\label{3.2}
A(x,t)= \F{1}{a(x)}
\int_{- \infty}^{x} \rho(x,t) - 2  \delta dx+ C(x)
\end{equation}
is a uniformly bounded function and satisfies
\begin{equation}\begin{array}{ll}
\label{3.3}
A_{t}(x,t)= \F{1}{a(x)} \int_{- \infty}^{x} \rho_{t}(x,t)dx
= \F{1}{a(x)} ( \E \rho_{x}- (\rho- 2 \delta)u),
\end{array}
\end{equation}
\begin{equation}\begin{array}{ll}
\label{3.4}
A_{x}(x,t)= -\F{a'(x)}{a^{2}(x)} \int_{- \infty}^{x} \rho(x,t) - 2  \delta dx
+ \F{1}{a(x)} ( \rho - 2  \delta)+ C'(x) \\\\
= B(x,t)+ \F{1}{a(x)} ( \rho - 2  \delta)+ C'(x) \geq 0,
\end{array}
\end{equation}
\begin{equation}
\label{3.5}
B(x,t)= -\F{a'(x)}{a^{2}(x)} \int_{- \infty}^{x} \rho(x,t) - 2  \delta dx \geq 0
\end{equation}
and
\begin{equation}\begin{array}{ll}
\label{3.6}
A_{xx}(x,t)= -(\F{a'(x)}{a^{2}(x)})' \int_{- \infty}^{x} \rho(x,t) - 2  \delta dx \\\\
-2 \F{a'(x)}{a^{2}(x)} ( \rho - 2  \delta)
+ \F{1}{a(x)} \rho_{x}+ C''(x).
\end{array}
\end{equation}
From (\ref{2.13}) and (\ref{2.14}), we obtain that
\begin{equation}\begin{array}{ll}
\label{3.7}
 z_{t}+ \lam^{\delta}_{1}z_{x}= \F{a(x)}{2 \tau} ((w-M-\tau A(x,t))-(z-M+ \tau A(x,t)) \\\\
 +  \E z_{xx}+ \F{2 \E}{\rho} \rho_{x}z_{x}- \F{\E}{2 \rho^{2}
   \sqrt{P'(\rho)}}(2P'+ \rho P'') \rho_{x}^{2}
   \end{array}
\end{equation}
and
\begin{equation}\begin{array}{ll}
\label{3.8}
 w_{t}+ \lam^{\delta}_{2}w_{x}= - \F{a(x)}{2 \tau} ((w-M-\tau A(x,t))-(z-M+ \tau A(x,t)) \\\\
+\E w_{xx}+ \F{2 \E}{\rho} \rho_{x}w_{x}- \F{\E}{2 \rho^{2}
   \sqrt{P'(\rho)}}(2P'+ \rho P'') \rho_{x}^{2},
\end{array}
\end{equation}
where $M>0$ is a suitable large constant.

Make the transformation
\begin{equation}
\label{3.9}
z=z_{1}+M - \tau A(x,t),  \quad  w=w_{1}+M +\tau A(x,t).
\end{equation}
Then the left-hand side of (\ref{3.7}) can be rewritten as
\begin{equation}\begin{array}{ll}
\label{3.10}
 L_{z_{1}}= z_{t}+ \lam^{\delta}_{1}z_{x}= z_{1t}- \tau A_{t}(x,t)+ \lam^{\delta}_{1}(z_{1x}- \tau A_{x}(x,t))  \\\\
= z_{1t}+ \lam^{\delta}_{1}z_{1x}- \F{ \tau}{a(x)} ( \E \rho_{x}- (\rho- 2 \delta)u) \\\\
  - \tau (u- \F{\rho-2 \delta}{\rho} \sqrt{P'( \rho)}) (B(x,t)+ \F{1}{a(x)} (\rho- 2 \delta)+C'(x))\\\\
  = z_{1t}+ \lam^{\delta}_{1}z_{1x} - \F{\E \tau}{a(x)} \rho_{x}
  + \F{\tau}{a(x)} \F{(\rho-2 \delta)^{2}}{\rho} \sqrt{P'( \rho)}  \\\\
  - \tau (B(x,t)+C'(x))( \int_{l}^{ \rho} \F{ \sqrt{P'(s)}}{s}ds- \F{\rho-2 \delta}{\rho} \sqrt{P'( \rho)}
  - (z_{1}- \tau A(x,t)+M)     ) \\\\
  = z_{1t}+ \lam^{\delta}_{1}z_{1x} + \tau (B(x,t)+C'(x)) z_{1} - \F{\E \tau}{a(x)} \rho_{x}
  + \F{\tau}{a(x)} \F{(\rho-2 \delta)^{2}}{\rho} \sqrt{P'( \rho)}  \\\\
  - \tau (B(x,t)+C'(x))( \int_{l}^{ \rho} \F{ \sqrt{P'(s)}}{s}ds- \F{\rho-2 \delta}{\rho} \sqrt{P'( \rho)})  \\\\
  + \tau (B(x,t)+C'(x)) (M- \tau A(x,t))
  \end{array}
\end{equation}
due to $ u= \int_{l}^{ \rho} \F{ \sqrt{P'(s)}}{s}ds-z$, and the following terms on the  right-hand side of (\ref{3.7})
\begin{equation}\begin{array}{ll}
\label{3.11}
 R_{z_{1}}= \E z_{xx} + \F{2 \E}{\rho} \rho_{x} z_{x}-
  \F{\E}{2 \rho^{2}  \sqrt{P'(\rho)}}(2P'+ \rho P'') \rho_{x}^{2}  \\\\
 =\E z_{1xx}- \E \tau A_{xx}+ \F{2 \E}{\rho} \rho_{x}( z_{1x}- \tau A_{x} )
 - \F{\E}{2 \rho^{2}  \sqrt{P'(\rho)}}(2P'+ \rho P'') \rho_{x}^{2} \\\\
 = \E z_{1xx} + \F{2 \E}{\rho} \rho_{x} z_{1x} - \F{\E}{2 \rho^{2}  \sqrt{P'(\rho)}}(2P'+ \rho P'') \rho_{x}^{2} - \F{2 \E \tau }{\rho} \rho_{x} A_{x}\\\\
 - \E \tau ( -(\F{a'(x)}{a^{2}(x)})' \int_{- \infty}^{x} \rho(x,t) - 2  \delta dx
-2 \F{a'(x)}{a^{2}(x)} ( \rho - 2  \delta)
+ \F{1}{a(x)} \rho_{x} )   \\\\
= \E z_{1xx} - \F{\E}{2 \rho^{2}
   \sqrt{P'(\rho)}}(2P'+ \rho P'') [\rho_{x}^{2} +  \F{4 \rho \sqrt{P'(\rho)}}{2P'+ \rho P''} \rho_{x}\tau A_{x} +( \F{2 \rho \sqrt{P'(\rho)}}{2P'+ \rho P''} \tau A_{x})^{2} ] \\\\
   + \F{2 \E}{\rho} \rho_{x} z_{1x} + \F{2 \E
\sqrt{P'(\rho)}}{2P'+ \rho P''} (\tau A_{x})^{2} +
\E \tau \F{a''(x)}{a^{2}(x)} \int_{- \infty}^{x} \rho(x,t) - 2  \delta dx \\\\
- 2 \E \tau \F{(a'(x))^{2}}{a^{3}(x)} \int_{- \infty}^{x} \rho(x,t) - 2  \delta dx
+ \E \tau (  2 \F{a'(x)}{a^{2}(x)} ( \rho - 2  \delta) -\F{1}{a(x)} \rho_{x}+C''(x) ) \\\\
\leq \E z_{1xx} + \F{2 \E}{\rho} \rho_{x} z_{1x} + \F{2 \E
\sqrt{P'(\rho)}}{2P'+ \rho P''} (\tau A_{x})^{2} \\\\+
\E \tau \F{a''(x)}{a^{2}(x)} \int_{- \infty}^{x} \rho(x,t) - 2  \delta dx
- \F{\E \tau}{a(x)} \rho_{x}+ \E \tau C''(x).
\end{array}
\end{equation}
When $ \gamma \geq 3$, we choose $l= 2 \delta$ in (\ref{1.11}). Since
\begin{equation}
\label{3.12}
\F{1}{\theta} ( \rho- 2 \delta) \F{\sqrt{P'(\rho)}}{\rho} \leq \int_{2 \delta}^{ \rho} \F{ \sqrt{P'(s)}}{s}ds \leq  ( \rho- 2 \delta) \F{\sqrt{P'(\rho)}}{\rho}
\quad \mbox{ for } \quad \gamma \geq 3,
\end{equation}
the following term in (\ref{3.10})
\begin{equation}
\label{3.13}
- \tau (B(x,t)+C'(x)) ( \int_{l}^{ \rho} \F{ \sqrt{P'(s)}}{s}ds- \F{\rho-2 \delta}{\rho} \sqrt{P'( \rho)})
\geq 0.
\end{equation}
By simple calculations,
\begin{equation}
\label{3.14}
A_{x}^{2} \leq 3 ( (-\F{a'(x)}{a^{2}(x)} \int_{- \infty}^{x} \rho(x,t) - 2  \delta dx)^{2}
+ \F{1}{a^{2}(x)} ( \rho - 2  \delta)^{2}+(C'(x))^{2})
\end{equation}
and
\begin{equation}
\label{3.15}
\F{\tau}{a(x)} \F{(\rho-2 \delta)^{2}}{\rho} \sqrt{P'( \rho)}  - \F{4 \E
\sqrt{P'(\rho)}}{2P'+ \rho P''} \tau^{2}  \F{1}{a^{2}(x)} ( \rho - 2  \delta)^{2} \geq 0,
\end{equation}
for $ \tau < 1$ and suitable small $\E$, in the range $ \rho \in ( 2 \delta, 2 M] $ (Note:  Finally we obtain the upper bound $M$ of $ \rho$ as expressed in (\ref{1.28})), and
\begin{equation}\begin{array}{ll}
\label{3.16}
\tau (B(x,t)+C'(x)) (M- \tau A(x,t))  -
\E \tau \F{a''(x)}{a^{2}(x)} \int_{- \infty}^{x} \rho(x,t) - 2  \delta dx \\\\
- \F{6 \E
\sqrt{P'(\rho)}}{2P'+ \rho P''} \tau^{2} ((-\F{a'(x)}{a^{2}(x)} \int_{- \infty}^{x} \rho(x,t) - 2  \delta dx)^{2}+ C'^{2}(x)) - \E \tau C''(x) \geq 0,
\end{array}
\end{equation}
where we choose $ \E$ to be much smaller than $\delta$ and assume
 \begin{equation}
\label{3.17}
\E |a''(x)| \leq |a'(x)|,  \quad \E |C''(x)| \leq C'(x).
\end{equation}
In fact, we may replace $a(x), C(x) $ by $ a^{\E}(x)= a(x)*G^{ \E}, C^{\E}(x)= C(x)*G^{ \E}$, $G^{ \E}$ being
a mollifier, so that  (\ref{3.17}) always true.

Then, from (\ref{3.14})-(\ref{3.17}), the following terms in (\ref{3.10}) and (\ref{3.11}) satisfy
\begin{equation}\begin{array}{ll}
\label{3.18}
D_{z_{1}}= \F{\tau}{a(x)} \F{(\rho-2 \delta)^{2}}{\rho} \sqrt{P'( \rho)}
  + \tau (B(x,t)+C'(x)) (M- \tau A(x,t))\\\\
  - \F{2 \E
\sqrt{P'(\rho)}}{2P'+ \rho P''} (\tau A_{x})^{2} -
\E \tau \F{a''(x)}{a^{2}(x)} \int_{- \infty}^{x} \rho(x,t) - 2  \delta dx - \E \tau C''(x) \geq 0.
\end{array}
\end{equation}
We obtain from (\ref{2.13}),(\ref{3.10}),(\ref{3.11}), (\ref{3.13}) and (\ref{3.18}) that
\begin{equation}\begin{array}{ll}
\label{3.19}
z_{1t}+ \lam^{\delta}_{1}z_{1x} + \tau (B(x,t)+C'(x)) z_{1}  \leq \F{a(x)}{ 2 \tau}(w_{1}-z_{1})
+ \E z_{1xx} + \F{2 \E}{\rho} \rho_{x} z_{1x}.
\end{array}
\end{equation}
When $ \gamma \in (1,3)$, we choose $l=0$ in (\ref{1.11}). The following term in (\ref{3.10})
\begin{equation}\begin{array}{ll}
\label{3.20}
- \tau (B(x,t)+C'(x))( \int_{l}^{ \rho} \F{ \sqrt{P'(s)}}{s}ds- \F{\rho-2 \delta}{\rho} \sqrt{P'( \rho)})
\geq - \tau (B(x,t)+C'(x)) \F{3- \gamma}{\gamma-1} \rho^{\theta} \\\\
= - \F{3- \gamma}{4} \tau (B(x,t)+C'(x)) (w_{1}+z_{1}+ 2 M).
\end{array}
\end{equation}
Since
\begin{equation}\begin{array}{ll}
\label{3.21}
D_{z_{1}} - \F{3- \gamma}{4} \tau (B(x,t)+C'(x))  2 M  \\\\
\geq \tau (B(x,t)+C'(x)) (M- \tau A(x,t)) - \F{3- \gamma}{4} \tau (B(x,t)+C'(x))  2 M \\\\
- \F{6 \E
\sqrt{P'(\rho)}}{2P'+ \rho P''} \tau^{2} ((-\F{a'(x)}{a^{2}(x)} \int_{- \infty}^{x} \rho(x,t) - 2  \delta dx)^{2}+ C'^{2}(x)) - \E \tau C''(x) \\\\
- \E \tau \F{a''(x)}{a^{2}(x)} \int_{- \infty}^{x} \rho(x,t) - 2  \delta dx \\\\
= \F{ \gamma-1}{2} \tau (B(x,t)+C'(x)) M- \tau^{2} (B(x,t)+C'(x))  A(x,t)  \\\\
- \F{6 \E
\sqrt{P'(\rho)}}{2P'+ \rho P''} \tau^{2} ((-\F{a'(x)}{a^{2}(x)} \int_{- \infty}^{x} \rho(x,t) - 2  \delta dx)^{2}+ C'^{2}(x)) - \E \tau C''(x) \\\\
- \E \tau \F{a''(x)}{a^{2}(x)} \int_{- \infty}^{x} \rho(x,t) - 2  \delta dx \geq 0,
\end{array}
\end{equation}
for a small $ \tau  $ or a sufficiently large $M$, then we
have
\begin{equation}\begin{array}{ll}
\label{3.22}
z_{1t}+ \lam^{\delta}_{1}z_{1x} + \tau (B(x,t)+C'(x)) z_{1} - \F{3- \gamma}{4} \tau (B(x,t)+C'(x)) (w_{1}+z_{1})
\\\\ \leq \F{a(x)}{ 2 \tau}(w_{1}-z_{1})
+ \E z_{1xx} + \F{2 \E}{\rho} \rho_{x} z_{1x}.
\end{array}
\end{equation}
Similarly, by using (\ref{2.14}), we have
\begin{equation}\begin{array}{ll}
\label{3.23}
w_{1t}+ \lam^{\delta}_{2}w_{1x} + \tau (B(x,t)+C'(x)) w_{1}  \leq \F{a(x)}{ 2 \tau }(z_{1}-w_{1})
+ \E w_{1xx} + \F{2 \E}{\rho} \rho_{x} w_{1x}
\end{array}
\end{equation}
when $ \gamma \geq 3$, and
\begin{equation}\begin{array}{ll}
\label{3.24}
w_{1t}+ \lam^{\delta}_{2}w_{1x} + \tau (B(x,t)+C'(x)) w_{1} - \F{3- \gamma}{4} \tau (B(x,t)+C'(x)) (w_{1}+z_{1})
\\\\ \leq \F{a(x)}{ 2 \tau}(z_{1}-w_{1})
+ \E w_{1xx} + \F{2 \E}{\rho} \rho_{x} w_{1x}
\end{array}
\end{equation}
when $ \gamma \in (1,3)$.

By applying the maximum principle to (\ref{3.19}) and (\ref{3.23}) when $ \gamma \geq 3$, and
to (\ref{3.22}) and (\ref{3.24}) when $ 1 < \gamma < 3$, we obtain $ z_{1} \leq 0, w_{1} \leq 0 $ or
\begin{equation}
\label{3.25}
z \leq M - \tau A(x,t) \leq M_{1},  \quad  w \leq M +\tau A(x,t) \leq M_{1},
\end{equation}
which reduce the estimates in (\ref{1.28}).

Finally, when $ \gamma=1$,  the perturbation pressure in (\ref{1.8}) is
\begin{equation}
\label{3.26}
 P_{1}( \rho, \delta)= \int^{\rho} \F{t-2
\delta}{t}P'(t)dt= \rho- 2 \delta  \ln \rho.
\end{equation}
Two eigenvalues of system (\ref{1.6}) are
\begin{equation}
\label{3.27} \lam^{\delta}_{1}= \F{m}{ \rho}- \F{\rho-2
\delta}{\rho} , \quad \lam^{\delta}_{2}= \F{m}{
\rho}+ \F{\rho-2 \delta}{\rho}
\end{equation}
with corresponding two same Riemann invariants
\begin{equation}
\label{3.28} z(\rho,u)=  \ln \rho -u,
\quad w(\rho,u)=\ln \rho +u.
\end{equation}
(\ref{3.10}) and (\ref{3.11}) are rewritten as follows:
\begin{equation}\begin{array}{ll}
\label{3.29}
 L_{z_{1}}
  = z_{1t}+ \lam^{\delta}_{1}z_{1x} + \tau (B(x,t)+C'(x)) z_{1} - \F{\E \tau}{a(x)} \rho_{x}
  + \F{\tau}{a(x)} \F{(\rho-2 \delta)^{2}}{\rho}  \\\\
  - \tau (B(x,t)+C'(x))( \ln \rho - \F{\rho-2 \delta}{\rho})
  + \tau (B(x,t)+C'(x)) (M- \tau A(x,t))
  \end{array}
\end{equation}
and
\begin{equation}\begin{array}{ll}
\label{3.30}
 R_{z_{1}}
\leq \E z_{1xx} + \F{2 \E}{\rho} \rho_{x} z_{1x} + \E (\tau A_{x})^{2} \\\\+
\E \tau \F{a''(x)}{a^{2}(x)} \int_{- \infty}^{x} \rho(x,t) - 2  \delta dx
- \F{\E \tau}{a(x)} \rho_{x}+ \E \tau C''(x),
\end{array}
\end{equation}
where $A_{x}$ satisfies (\ref{3.14}). In this case, (\ref{3.15})  is changed to
\begin{equation}
\label{3.31}
\F{\tau}{a(x)} \F{(\rho-2 \delta)^{2}}{\rho}   - 2 \E
 \tau^{2}  \F{1}{a^{2}(x)} ( \rho - 2  \delta)^{2} \geq 0,
\end{equation}
in the range $ \rho \in ( 2 \delta, 2 M]$, for small $\E$.

Now we separate the points $(x,t)$ into two different parts. First, at the points $(x,t)$, where $ \rho(x,t) \leq 1 $,
we have
\begin{equation}
\label{3.32}
- \tau B(x,t)( \ln \rho - \F{\rho-2 \delta}{\rho}) \geq 0.
\end{equation}
Then (\ref{3.21}) is changed to
\begin{equation}\begin{array}{ll}
\label{3.33}
D_{z_{1}}= \F{\tau}{a(x)} \F{(\rho-2 \delta)^{2}}{\rho}
  + \tau (B(x,t)+C'(x)) (M- \tau A(x,t))\\\\
  -  \E (\tau A_{x})^{2} -
\E \tau \F{a''(x)}{a^{2}(x)} \int_{- \infty}^{x} \rho(x,t) - 2  \delta dx - \E \tau C''(x) \geq 0,
\end{array}
\end{equation}
and we obtain
\begin{equation}\left\{
\label{3.34}
\begin{array}{l}
z_{1t}+ \lam^{\delta}_{1}z_{1x} + \tau (B(x,t)+C'(x)) z_{1}  \leq \F{a(x)(\rho- 2 \delta)}{ 2 \tau \rho}(w_{1}-z_{1})
+ \E z_{1xx} + \F{2 \E}{\rho} \rho_{x} z_{1x},   \\\\
w_{1t}+ \lam^{\delta}_{2}w_{1x} + \tau (B(x,t)+C'(x)) w_{1}  \leq \F{a(x)(\rho- 2 \delta)}{ 2 \tau \rho}(z_{1}-w_{1})
+ \E w_{1xx} + \F{2 \E}{\rho} \rho_{x} w_{1x}.
\end{array}\right.
\end{equation}
Second, at the points $(x,t)$, where $ \rho(x,t) > 1 $,
\begin{equation}\begin{array}{ll}
\label{3.35}
- \tau (B(x,t)+C'(x)) \ln \rho
= - \F{1}{2} \tau (B(x,t)+C'(x)) (w_{1}+z_{1}+ 2 M) \\\\
= - \F{1}{2} \tau (B(x,t)+C'(x)) (w_{1}+z_{1}) - \tau M (B(x,t)+C'(x)),
\end{array}
\end{equation}
\begin{equation}
\label{3.36}
 \tau (B(x,t)+C'(x)) \F{\rho-2 \delta}{\rho} \geq \F{1}{2} \tau (B(x,t)+C'(x)).
\end{equation}
Then
\begin{equation}\begin{array}{ll}
\label{3.37}
D_{z_{1}}- \tau M (B(x,t)+C'(x)) + \F{1}{2} \tau (B(x,t)+C'(x))  \\\\
 \geq \F{1}{2} \tau (B(x,t)+C'(x)) -
   \tau^{2} (B(x,t)+C'(x))  A(x,t)  \\\\
  - 3 \E
\tau^{2} ((-\F{a'(x)}{a^{2}(x)} \int_{- \infty}^{x} \rho(x,t) - 2  \delta dx)^{2}+ C'^{2}(x))
- \E \tau C''(x) \geq 0
\end{array}
\end{equation}
and
\begin{equation}\left\{
\label{3.38}
\begin{array}{l}
z_{1t}+ \lam^{\delta}_{1}z_{1x} + \tau (B(x,t)+C'(x)) z_{1} - \F{1}{2} \tau (B(x,t)+C'(x)) (w_{1}+z_{1})
\\\\ \leq \F{a(x)(\rho- 2 \delta)}{ 2 \tau \rho}(w_{1}-z_{1})
+ \E z_{1xx} + \F{2 \E}{\rho} \rho_{x} z_{1x},   \\\\
w_{1t}+ \lam^{\delta}_{2}w_{1x} + \tau (B(x,t)+C'(x)) w_{1} - \F{1}{2} \tau (B(x,t)+C'(x)) (w_{1}+z_{1})
\\\\ \leq \F{a(x)(\rho- 2 \delta)}{ 2 \tau \rho}(z_{1}-w_{1})
+ \E w_{1xx} + \F{2 \E}{\rho} \rho_{x} w_{1x}.
\end{array}\right.
\end{equation}
Uniting (\ref{3.34}) and (\ref{3.38}) together, we obtain the following two inequalities at any point $(x,t) \in (-\infty,\infty) \times (0, \infty)$,
\begin{equation}\left\{
\label{3.39}
\begin{array}{l}
z_{1t}+ \lam^{\delta}_{1}z_{1x} + \tau (B(x,t)+C'(x)) z_{1} + l_{0}(x,t) (w_{1}+z_{1})
\\\\ \leq \F{a(x)(\rho- 2 \delta)}{ 2 \tau \rho}(w_{1}-z_{1})
+ \E z_{1xx} + \F{2 \E}{\rho} \rho_{x} z_{1x},   \\\\
w_{1t}+ \lam^{\delta}_{2}w_{1x} + \tau (B(x,t)+C'(x)) w_{1} +l_{0}(x,t) (w_{1}+z_{1})
\\\\ \leq \F{a(x)(\rho- 2 \delta)}{ 2 \tau \rho}(z_{1}-w_{1})
+ \E w_{1xx} + \F{2 \E}{\rho} \rho_{x} w_{1x},
\end{array}\right.
\end{equation}
where $ l_{0}(x,t) \leq 0$ is a suitable function.

By applying the maximum principle to (\ref{3.39}), we obtain $ z_{1} \leq 0, w_{1} \leq 0 $ or
\begin{equation}
\label{3.40}
\ln \rho^{\E,\delta}- u^{\E,\delta}  \leq M - \tau A(x,t),  \quad  \ln \rho^{\E,\delta}+ u^{\E,\delta} \leq M +\tau A(x,t),
\end{equation}
which reduce the estimates in (\ref{1.29}).

Thus we complete the proof of Theorem 2.

\section{ \bf Proof of Theorem 3}

To study the relaxation limit, namely the limit of  $(N^{\tau}(x,s),J^{\tau}(x,s),\Up^{\tau}(x,s))$ as $ \tau \rightarrow 0^{+}, \E\rightarrow 0^{+}, \delta \rightarrow 0^{+}$, we add also a small perturbation $\delta$  to the terms $ \rho E- \F{1}{\tau} a(x) \rho u $ in the second equation in (\ref{1.1}) or replace (\ref{1.6}) by the following system
\begin{equation}\left\{
\label{4.1}
\begin{array}{l}
 \rho_{t}+( (\rho-2 \delta) u)_{x} =\E \rho_{xx},         \\\\
 ( \rho u)_t+( \rho u^2- \delta u^{2}+ P_{1}( \rho, \delta))_x= \E
 (\rho u)_{xx}+ (\rho- 2 \delta) E- \F{1}{\tau} a(x) (\rho- 2 \delta) u, \\\\
 E_x= (\rho- 2 \delta) -b(x).
  \end{array}\right.
\end{equation}
Repeating the proof given in the last section, we may prove that the solutions $(\rho^{\E,\delta},m^{\E,\delta},E^{\E,\delta})$,
of the problem (\ref{4.1}),(\ref{1.3}) and (\ref{1.7}), satisfy
\begin{equation}\left\{
\label{4.2}
\begin{array}{l}
z_{1t}+ \lam^{\delta}_{1}z_{1x} + \tau (B(x,t)+C'(x)) z_{1}  \leq \F{a(x)(\rho- 2 \delta)}{ 2 \tau \rho}(w_{1}-z_{1})
+ \E z_{1xx} + \F{2 \E}{\rho} \rho_{x} z_{1x},   \\\\
w_{1t}+ \lam^{\delta}_{2}w_{1x} + \tau (B(x,t)+C'(x)) w_{1}  \leq \F{a(x)(\rho- 2 \delta)}{ 2 \tau \rho}(z_{1}-w_{1})
+ \E w_{1xx} + \F{2 \E}{\rho} \rho_{x} w_{1x}
\end{array}\right.
\end{equation}
when $ \gamma \geq 3$;
\begin{equation}\left\{
\label{4.3}
\begin{array}{l}
z_{1t}+ \lam^{\delta}_{1}z_{1x} + \tau (B(x,t)+C'(x)) z_{1} - \F{3- \gamma}{4} \tau (B(x,t)+C'(x)) (w_{1}+z_{1})
\\\\ \leq \F{a(x)(\rho- 2 \delta)}{ 2 \tau \rho}(w_{1}-z_{1})
+ \E z_{1xx} + \F{2 \E}{\rho} \rho_{x} z_{1x},   \\\\
w_{1t}+ \lam^{\delta}_{2}w_{1x} + \tau (B(x,t)+C'(x)) w_{1} - \F{3- \gamma}{4} \tau (B(x,t)+C'(x)) (w_{1}+z_{1})
\\\\ \leq \F{a(x)(\rho- 2 \delta)}{ 2 \tau \rho}(z_{1}-w_{1})
+ \E w_{1xx} + \F{2 \E}{\rho} \rho_{x} w_{1x}
\end{array}\right.
\end{equation}
when $ \gamma \in (1,3)$ and
\begin{equation}\left\{
\label{4.4}
\begin{array}{l}
z_{1t}+ \lam^{\delta}_{1}z_{1x} + \tau (B(x,t)+C'(x)) z_{1} +l_{0}(x,t) (w_{1}+z_{1})
\\\\ \leq \F{a(x)(\rho- 2 \delta)}{ 2 \tau \rho}(w_{1}-z_{1})
+ \E z_{1xx} + \F{2 \E}{\rho} \rho_{x} z_{1x},   \\\\
w_{1t}+ \lam^{\delta}_{2}w_{1x} + \tau (B(x,t)+C'(x)) w_{1} +l_{0}(x,t) (w_{1}+z_{1})
\\\\ \leq \F{a(x)(\rho- 2 \delta)}{ 2 \tau \rho}(z_{1}-w_{1})
+ \E w_{1xx} + \F{2 \E}{\rho} \rho_{x} w_{1x}
\end{array}\right.
\end{equation}
when $\gamma=1$, where $l_{0}(x,t) \leq 0$ and $w,z$ are given in (\ref{1.11}).

Then we may obtain the upper bound
\begin{equation}
\label{4.5}
  w(\rho^{\E, \delta},u^{\E, \delta}) \leq M, \quad z(\rho^{\E, \delta},u^{\E, \delta}) \leq M,
  \end{equation}
where $M$ is a suitable large positive constant, which depend only on the bound of the initial data, but is independent of $\E, \delta, \tau$ and the time $t$.

When $P(\rho)$ satisfies the condition (\ref{1.14}) in Theorem 1, which is corresponding to the case of $\gamma > 1$, we have
\begin{equation}
\label{4.6}
  0 < 2 \delta \leq \rho^{\E, \delta} \leq M_{1}, \quad |u^{\E, \delta}| \leq M_{1}
  \end{equation}
 or when $\gamma=1$,
 \begin{equation}\left\{
\label{4.7}
\begin{array}{l}
0 < 2 \delta \leq \rho^{\E,\delta} \leq M, \quad \ln \rho^{\E,\delta}-u^{\E,\delta} \leq M,
\quad \ln \rho^{\E,\delta}+u^{\E,\delta} \leq M,    \\\\
  -M_{1} \leq  \rho^{\E,\delta} (\ln \rho^{\E,\delta}- M)
\leq \rho^{\E,\delta}u^{\E,\delta} \leq \rho^{\E,\delta}(M- \ln \rho^{\E,\delta}) \leq M_{1},  \\\\
 0 \leq  \rho^{\E,\delta} (u^{\E,\delta})^{2}
\leq  \rho^{\E,\delta} \max \{(\ln \rho^{\E,\delta}- M)^{2}, (\ln \rho^{\E,\delta}+ M)^{2} \} \leq M_{2}
  \end{array}\right.
\end{equation}
for two suitable positive constants $M_{1}$ and $M_{2}$.

As did in \cite{MN2}, we introduce the scaled variables in  (\ref{4.1})
\begin{equation}
\label{4.8}
N^{\tau}(x,s)= \rho^{\E,\delta}(x,\F{s}{\tau}),  \quad J^{\tau}(x,s)=\F{1}{\tau} m^{\E,\delta}(x,\F{s}{\tau}),
\quad \Up^{\tau}(x,s)= E^{\E,\delta}(x,\F{s}{\tau}),
  \end{equation}
then (\ref{4.1}) is rewritten as
\begin{equation}\left\{
\label{4.9}
\begin{array}{l}
 N^{\tau}_{s}+( (N^{\tau}-2 \delta) U^{\tau})_{x}
 =\F{\E}{\tau} N^{\tau}_{xx},         \\\\
\tau^{2} J^{\tau}_{s}+(\tau^{2} (N^{\tau} (U^{\tau})^2
- \delta (U^{\tau})^2)+ P_{1}( N^{\tau}, \delta))_x\\\\= \E \tau
 J^{\tau}_{xx}+(N^{\tau}- 2 \delta) \Up^{\tau}- a(x) (N^{\tau}- 2 \delta) U^{\tau}, \\\\
\Up^{\tau}_x= (N^{\tau}- 2 \delta) -b(x),
  \end{array}\right.
\end{equation}
where $ U^{\tau}(x,s)= \F{J^{\tau}(x,s)}{N^{\tau}(x,s)} = \F{1}{\tau} u^{\E,\delta}(x,\F{s}{\tau})$.

Let
\begin{equation}
\label{4.10}
      (\eta^{\star}(\rho,m),\ q^{\star}(\rho,m))=\Big(\F{m^2}{2 \rho} + \rho \int_{2 \delta}^{\rho} \F{P(s)}{s^2}ds,\
\F{m^3}{2 \rho^2}+( \F{P( \rho)}{\rho}+ \int_{2 \delta}^{\rho} \F{P(s)}{s^2}ds)m\Big).
\end{equation}
Our approximated solutions satisfy (\ref{4.1}) or equivalently
\begin{equation}\left\{
\label{4.11}
\begin{array}{l}
\rho_{t}+(\rho u)_{x}-2 \delta u_{x}=  \E \rho_{xx},     \\  \\
       ( \rho u)_t+( \rho u^2+ P( \rho))_x - 2 \delta u u_{x} - 2 \delta \F{P'(\rho)}{\rho} \rho_{x}\\\\
      = \E (\rho u)_{xx} + (\rho- 2 \delta) E- a(x) \F{(\rho- 2 \delta) u}{\tau}.
\end{array}\right.
\end{equation}
Then
\begin{equation}\begin{array}{ll}
\label{4.12}
\eta^{\star}_{t}( \rho^{\E, \delta},u^{\E, \delta})+ q^{\star}_{x}( \rho^{\E, \delta},u^{\E, \delta})
- \F{1}{3} \delta ((u^{\E, \delta})^{3})_{x} - 2 \delta (( \int_{2 \delta}^{\rho^{\E, \delta}} \F{P(s)}{s^2}ds+ \F{P(\rho^{\E, \delta})}{\rho^{\E, \delta}}) u^{\E, \delta})_{x} \\\\
= \E \eta^{\star}_{xx}( \rho^{\E, \delta},u^{\E, \delta})  + (\rho^{\E, \delta}- 2 \delta) u^{\E, \delta} E^{\E, \delta}
- a(x) \F{(\rho^{\E, \delta}- 2 \delta) (u^{\E, \delta})^{2}}{\tau} \\\\
-\E ( \rho^{\E, \delta}_{x}, m^{\E, \delta}_{x}) \cdot \nabla^{2}
 \eta^{\star}( \rho^{\E, \delta}, m^{\E, \delta}) \cdot ( \rho^{\E, \delta}_{x}, m^{\E, \delta}_{x})^{T}.
\end{array}
\end{equation}
Let
\begin{equation}
\label{4.13}
\eta(N^{\tau},U^{\tau})
=  \F{1}{2} \tau^{2} N^{\tau}  (U^{\tau})^{2}
 + N^{\tau} \int_{2 \delta}^{ N^{\tau} } \F{P(s)}{s^{2}} ds.
\end{equation}
Then we have from  (\ref{4.12}) that
\begin{equation}\begin{array}{ll}
\label{4.14}
 \eta_{s}(N^{\tau},U^{\tau})+ \F{1}{\tau} Q_{x}(N^{\tau},U^{\tau})
= - \F{\E}{\tau} ( \rho^{\E}_{x}, m^{\E}_{x}) \cdot \nabla^{2}
 \eta^{\star}( \rho^{\E}, m^{\E}) \cdot ( \rho^{\E}_{x}, m^{\E}_{x})^{T} \\\\
+ \F{\E}{\tau} \eta_{xx}(N^{\tau},U^{\tau})
 + \Up^{\tau} (N^{\tau}- 2 \delta) U^{\tau} - a(x)  (N^{\tau}- 2 \delta)  (U^{\tau})^{2},
 \end{array}
\end{equation}
where
\begin{equation}
\label{4.15}
Q(N^{\tau},U^{\tau})=
q^{\star}( \rho^{\E},u^{\E})
- \F{1}{3} \delta (u^{\E})^{3} - 2 \delta ( \int_{2 \delta}^{\rho^{\E}} \F{P(s)}{s^2}ds+ \F{P(\rho^{\E})}{\rho^{\E}}) u^{\E}.
\end{equation}
Integrating (\ref{4.14}) in $R$, we have from (\ref{1.31}), (\ref{2.1}) and the boundedness of $ \rho^{\E}- 2 \delta$ in $L^{1}(R)$  that
\begin{equation}\begin{array}{ll}
\label{4.16}
\int_{-\infty}^{\infty}  \eta_{s}(N^{\tau},U^{\tau}) dx
+ a_{0} \int_{-\infty}^{\infty}  (N^{\tau}- 2 \delta)  (U^{\tau})^{2} dx \\\\
\leq
 \int_{-\infty}^{\infty} |\Up^{\tau} (N^{\tau}- 2 \delta) U^{\tau}| dx
\leq  M ( \int_{-\infty}^{\infty}  (N^{\tau}- 2 \delta)  (U^{\tau})^{2} dx )^{\F{1}{2}},
\end{array}
\end{equation}
where $a_{0}$ is the positive lower bound of $ a(x)$.

Letting $\phi(s)= \int_{-\infty}^{\infty} (N^{\tau}- 2 \delta)  (U^{\tau})^{2} dx$,
  we have from (\ref{4.16}) that
 \begin{equation}
\label{4.17}
\int_{0}^{L} \phi(s) ds \leq M_{1} \int_{0}^{L} \phi^{\F{1}{2}}(s) ds + M_{2}
\leq M_{1} L^{\F{1}{2}} (\int_{0}^{L} \phi(s) ds)^{\F{1}{2}} + M_{2},
\end{equation}
for two suitable constants $M_{1},M_{2}$. Then $ \int_{0}^{L} \phi(s) ds \leq M(L)$ or
\begin{equation}
\label{4.18}
 \int_{0}^{L} \int_{-\infty}^{\infty} (N^{\tau}- 2 \delta)  (U^{\tau})^{2} dx ds \leq M(L),
\end{equation}
where $ M(L)$ depends on $L$, but is independent of $\E,\delta,\tau$.

By simple calculations,
\begin{equation}\begin{array}{ll}
\label{4.19}
\F{\E}{\tau} ( \rho_{x}, m_{x}) \cdot \nabla^{2}
 \eta^{\star}( \rho, m) \cdot ( \rho_{x}, m_{x})^{T} \\\\
 = \F{\E}{\tau} [( \F{m^{2}}{\rho^{3}}+ \F{P'(\rho)}{\rho}) \rho_{x}^{2} - 2 \F{m}{\rho^{2}} \rho_{x} m_{x}
 + \F{1}{\rho} m_{x}^{2}]  \\\\
 = \F{\E}{\tau} [ \F{m^{2}}{\rho^{3}} \rho_{x}^{2} - 2 \F{m}{\rho^{2}} \rho_{x} m_{x}
 + \F{1}{\rho} m_{x}^{2}] + \F{\E}{\tau} \F{P'(\rho)}{\rho} \rho_{x}^{2}  \geq \F{\E}{\tau} \F{P'(\rho)}{\rho} \rho_{x}^{2}
\end{array}
\end{equation}
and also
\begin{equation}
  \begin{array}{ll}
  \label{4.20}
\F{\E}{\tau} ( \rho_{x}, m_{x}) \cdot \nabla^{2}
 \eta^{\star}( \rho, m) \cdot ( \rho_{x}, m_{x})^{T} \\\\
 =\F{\E}{\tau} [( \F{m^{2}}{\rho^{3}}+ \F{P'(\rho)}{\rho}) \rho_{x}^{2} - 2 \F{m}{\rho^{2}} \rho_{x} m_{x}
 + \F{1}{\F{m^{2}}{\rho^{3}}+ \F{P'(\rho)}{\rho}}(\F{m}{\rho^{2}})^{2} m_{x}^{2}] \\\\
 + \F{\E}{\tau} (\F{1}{\rho}- \F{1}{\F{m^{2}}{\rho^{3}}+ \F{P'(\rho)}{\rho}}(\F{m}{\rho^{2}})^{2})  m_{x}^{2}]
 \geq \F{\E}{\tau} (\F{1}{\rho}- \F{1}{\F{m^{2}}{\rho^{3}}+ \F{P'(\rho)}{\rho}}(\F{m}{\rho^{2}})^{2})  m_{x}^{2},
\end{array}
\end{equation}
where
\begin{equation}
\label{4.21}
\F{1}{\rho}- \F{1}{\F{m^{2}}{\rho^{3}}+ \F{P'(\rho)}{\rho}}(\F{m}{\rho^{2}})^{2}
= \F{1}{\rho} \F{P'(\rho)}{u^{2}+P'(\rho)} \geq c_{1} P'(2 \delta)
\end{equation}
 for a suitable constant $ c_{1} > 0$, and for general pressure $P(\rho)$; or
 \begin{equation}
\label{4.22}
\F{1}{\rho}- \F{1}{\F{m^{2}}{\rho^{3}}+ \F{P'(\rho)}{\rho}}(\F{m}{\rho^{2}})^{2} =
\F{1}{\rho}- \F{1}{\F{m^{2}}{\rho^{3}}+ \F{1}{\rho}}(\F{m}{\rho^{2}})^{2}
= \F{1}{\rho u^{2}+ \rho} \geq c_{1}
\end{equation}
when $P(\rho)= \rho$.

Then we have from (\ref{4.14}),(\ref{4.16}), (\ref{4.18})-(\ref{4.22}) that
\begin{equation}
\label{4.23}
\F{\E}{\tau} P'(2 \delta) (|(N^{\tau}_{x})^{2}|_{L^{1}(R \times [0,L]}+
  (J^{\tau}_{x})^{2}|_{L^{1}(R \times [0,L]}) \leq M(L),
\end{equation}
where $M(L)$ depends only on $L$.

Now we prove the following Lemma:
\begin{lemma}
The sequence of functions $( N^{\tau}(x,s),J^{\tau}(x,s), \Up^{\tau}(x,s))$ satisfies
\begin{equation}
\label{4.24}
 N^{\tau} \in L^{2}_{loc}(R \times R^{+}),
 \quad  (N^{\tau}-2 \delta) U^{\tau}  \in L^{2}_{loc}(R \times R^{+}),
 \end{equation}
 \begin{equation}
\label{4.25}
 \tau^{2} (N^{\tau} (U^{\tau})^2
- \delta (U^{\tau})^2)+ P_{1}( N^{\tau}, \delta) \in L^{2}_{loc}(R \times R^{+}),
  \end{equation}
and
\begin{equation}
\label{4.26}
 N^{\tau}_{s}+( (N^{\tau}-2 \delta) U^{\tau})_{x}, \quad \tau^{2} J^{\tau}_{s}+(\tau^{2} (N^{\tau} (U^{\tau})^2
- \delta (U^{\tau})^2)+ P_{1}( N^{\tau}, \delta))_x
\end{equation}
 are compact in $ H^{-1}_{loc}(R \times R^{+})$.
\end{lemma}
{\bf Proof of Lemma 5.} By applying (\ref{2.11}), (\ref{4.6}), (\ref{4.7}) and (\ref{4.18}), we have (\ref{4.24}) and (\ref{4.25})
immediately.

To prove (\ref{4.26}), first, for any $\varphi\in H^{1}_{0}(R \times R^{+})$, we have from (\ref{4.23}) and the first equation in (\ref{4.9}) that
\begin{equation}\begin{array}{ll}
\label{4.27}
|\int_{0}^{\infty} \int_{- \infty}^{\infty} \F{\E}{\tau} N^{\tau}_{xx} \varphi dx ds|
= |\int_{0}^{\infty} \int_{- \infty}^{\infty} \F{\E}{\tau \sqrt{P'( 2 \delta)}}  \sqrt{P'( 2 \delta)} N^{\tau}_{x} \varphi_{x} dx ds|
\rightarrow 0
\end{array}
\end{equation}
if we choose $ \E = o( \sqrt{P'( 2 \delta)} \tau) $ and $ \delta, \tau$ go to zero,
where $ \Omega $ is the compact support set of $ \varphi$. Then $ N^{\tau}_{s}+( (N^{\tau}-2 \delta) U^{\tau})_{x}$
 are compact in $ H^{-1}_{loc}(R \times R^{+})$.

 Second, since
 \begin{equation}\begin{array}{ll}
\label{4.28}
|\int_{0}^{\infty} \int_{- \infty}^{\infty} \E \tau J^{\tau}_{xx} \varphi dx ds|
= |\int_{0}^{\infty} \int_{- \infty}^{\infty} \F{\E \tau}{\sqrt{P'( 2 \delta)}} \sqrt{P'( 2 \delta)} J^{\tau}_{x} \varphi_{x} dx ds|
\rightarrow 0
\end{array}
\end{equation}
if we choose $ \E = o( \sqrt{P'( 2 \delta)} ) $ and $ \delta, \tau$ go to zero,
then $ \E \tau J^{\tau}_{xx} $  are compact in $ H^{-1}_{loc}(R \times R^{+})$.

Third, by using (\ref{4.18})
\begin{equation}\begin{array}{ll}
\label{4.29}
\int \int_{\Omega}
a(x) | (N^{\tau}- 2 \delta) U^{\tau}| dx ds  \\\\
\leq  M  \int \int_{\Omega}  |N^{\tau}- 2 \delta| dx ds )^{\F{1}{2}} \int \int_{\Omega}  (N^{\tau}- 2 \delta)  (U^{\tau})^{2} dx ds )^{\F{1}{2}} \leq M_{1},
\end{array}
\end{equation}
then the terms $ (N^{\tau}- 2 \delta) \Up^{\tau}- a(x) (N^{\tau}- 2 \delta) U^{\tau} $,
on the right-hand side of the second equation in (\ref{4.9}), are bounded in $L^{1}_{loc}(R \times R^{+})$, and so compact in $W^{-1,q}_{loc}(R \times R^{+}), q \in (1,2)$ by using the Sobolev compact embedding theorem.
Uniting (\ref{4.28}) and (\ref{4.29}), we have that the right-hand side of the second equation in (\ref{4.9}) is compact in $W^{-1,q}_{loc}(R \times R^{+}), q \in (1,2)$.
On the other side, the left-hand side, $ \tau^{2} J^{\tau}_{s}+(\tau^{2} (N^{\tau} (U^{\tau})^2
- \delta (U^{\tau})^2)+ P_{1}( N^{\tau}, \delta))_x $  is bounded in $W^{-1,\infty}(R \times R^{+})$.

Therefore, by using the  Murat embedding theorem \cite{Mu}, $ \tau^{2} J^{\tau}_{s}+(\tau^{2} (N^{\tau} (U^{\tau})^2
- \delta (U^{\tau})^2)+ P_{1}( N^{\tau}, \delta))_x $  are compact in $H^{-1}_{loc}(R \times R^{+})$.
{\bf Lemma 5 is proved}.

Let
\begin{equation}\begin{array}{ll}
\label{4.30}
(N^{\tau}, (N^{\tau}-2 \delta) U^{\tau}, \tau^{2} J^{\tau},\tau^{2} (N^{\tau} (U^{\tau})^2
- \delta (U^{\tau})^2)+ P_{1}( N^{\tau}, \delta)) \\\\
\rightharpoonup
(N,J,v_3,v_4), \quad \mbox{ weakly in } \quad L^2(\Omega).
\end{array}
\end{equation}
By using the Div-Curl lemma (see \cite{Mu,Ta}), we have from (\ref{4.24})-(\ref{4.26}) in Lemma 5 that
\begin{equation}\begin{array}{ll}
\label{4.31}
 N^{\tau} \times (\tau^{2} (N^{\tau} (U^{\tau})^2
- \delta (U^{\tau})^2)+ P_{1}( N^{\tau}, \delta))-
\tau^{2} J^{\tau}  \times (N^{\tau}-2 \delta) U^{\tau}  \\\\
= N^{\tau} P_{1}( N^{\tau}, \delta) + \delta \tau^{2} N^{\tau} (U^{\tau})^2
\rightharpoonup
N v_{4}-J v_{3}
\end{array}
\end{equation}
in the sense of distributions.

By using (\ref{4.6})-(\ref{4.7}),  $ \tau J^{\tau}, \tau^{2} N^{\tau} (U^{\tau})^2 $ are uniformly bounded. Then
\begin{equation}
\label{4.32}
  \delta \tau^{2} N^{\tau} (U^{\tau})^2,  \quad \tau^{2} J^{\tau} \rightarrow 0, \quad  a.e., \mbox{ as } \quad
 \tau,   \delta \rightarrow 0.
\end{equation}
Thus $ v_3=0$ and
\begin{equation}
\label{4.33}
N^{\tau} P_{1}( N^{\tau}, \delta)
\rightharpoonup
N v_{4}
\end{equation}
due to (\ref{4.31}).

Furthermore, by using (\ref{4.18}),
\begin{equation}
\begin{array}{ll}
\label{4.34}
|\tau^{2} (N^{\tau} (U^{\tau})^2
- \delta (U^{\tau})^2)|_{L^{2}_{loc}(R \times R^{+})}=
|\tau^{4} (N^{\tau}-  \delta)^{2} (U^{\tau})^4|_{L^{1}_{loc}(R \times R^{+})} \\\\
 \leq
2 |\tau^{4} (N^{\tau}- 2 \delta)^{2} (U^{\tau})^4|_{L^{1}_{loc}(R \times R^{+})}+
2 |\tau^{4}  \delta^{2} (U^{\tau})^4|_{L^{1}_{loc}(R \times R^{+})} \\\\
\leq M \tau^{2} | (N^{\tau}- 2 \delta) (U^{\tau})^2|_{L^{1}_{loc}(R \times R^{+})}+
2 |\tau^{4}  \delta^{2} (U^{\tau})^4|_{L^{1}_{loc}(R \times R^{+})}  \rightarrow 0,
\end{array}
\end{equation}
as $ \delta, \tau$ go to zero,
then
\begin{equation}
\label{4.35}
P_{1}( N^{\tau}, \delta)
\rightharpoonup
v_{4}, \quad \mbox{ weakly in } L^2(\Omega).
\end{equation}
Moreover, since for fixed $ \tau$,
 \begin{equation}
\label{4.36}
  P_{1}( N^{\tau}, \delta)
\rightarrow  P( N^{\tau}), \quad a.e.,  \mbox{ as }  \delta \rightarrow 0,
\end{equation}
 we have from (\ref{4.33}), (\ref{4.35}) and (\ref{4.36}) that
\begin{equation}
\label{4.37}
\overline{ N^{\tau} P( N^{\tau})}= N \overline{ P( N^{\tau}) } ,
\end{equation}
where $ \overline{f(u^{\tau})}$ denotes the weak-star limit of
$f(u^{\tau})$.

When $P(\rho)= \rho$, (\ref{4.37}) reduces the pointwise convergence of $N^{\tau}$ immediately. When $P(\rho)$ is a convex function, the well-known argument of Minty  shows that $ \overline{ P( N^{\tau}) } =P(N)$ (see \cite{MN2,MM}) and
\begin{equation}
\label{4.38}
N^{\tau} \rightarrow N,  \mbox{ strongly in }  L^{p}(\Omega),  \mbox{ for all }  p \geq 1.
\end{equation}
Let the weak star limit of $ \Up^{\tau}(x,s) $ be $\Up(x,s)$. By using the third equation in (\ref{4.9}), we
know that $ \Up^{\tau}_{x} $ are compact in $H^{-1}_{loc}(R \times R^{+})$. Moreover,
since
\begin{equation}
\label{4.39}
\Up^{\tau}_{s}= \int_{- \infty}^{x} N^{\tau}_{s} dx =  -(N^{\tau}-2 \delta) U^{\tau}
 + \F{\E}{\tau} N^{\tau}_{x}
\end{equation}
are bounded in $L^{1}_{loc}(R \times R^{+})$, and so compact in $W^{-1,q}_{loc}(R \times R^{+}), q \in (1,2)$,
then $ \Up^{\tau}_{s} $ are also compact in $H^{-1}_{loc}(R \times R^{+})$.

We may apply the Div-Curl lemma to the pairs of functions
\begin{equation}
\label{4.40}
(0,\Up^{\tau}),    \quad (\Up^{\tau},0),
\end{equation}
to obtain
\begin{equation}
 \label{4.41}
\overline{ \Up^{\tau}} \cdot \overline{ \Up^{\tau}}=  \overline{
(\Up^{\tau})^{2}},
\end{equation}
 which deduces the pointwise convergence of $\Up^{\tau}(x,s)$,
\begin{equation}
\label{4.42}
\Up^{\tau}(x,s) \rightarrow \Up(x,s),  \mbox{ strongly in }  L^{p}(\Omega),  \mbox{ for all }  p \geq 1.
\end{equation}
Finally, since $ \Up^{\tau}_{x} $ is uniform bounded,
then $ \Up(x,s) $ is Lipschitz continuous with respect to the space variable $x$. Letting $ \E, \delta, \tau$ in (\ref{4.9}) go to zero,
we have that the limit $(N,J,\Up)$ satisfies the  drift-diffusion equations (\ref{1.22}) in the sense of distributions.
Thus {\bf Theorem 3 is proved.}






\end{document}